\documentclass[preprint,11pt]{elsarticle}
\usepackage{amsmath}
\usepackage{amssymb}
\usepackage{latexsym}
\usepackage{graphicx}
\usepackage[dvips]{epsfig}
\usepackage[usenames]{color}
\usepackage{psfrag}

\input amssym.def
\newsymbol\rtimes 226F
\newfont{\nset}{msbm10}

\newtheorem{theo}{Theorem}[section]
\newtheorem{theorem}[theo]{Theorem}

\journal{Physica A}

\begin{document}

\begin{frontmatter}

\title{The Number of Spanning Trees of an  Infinite Family of Outerplanar, Small-World and Self-Similar Graphs}

\author{Francesc Comellas, Al\'icia Miralles}
\address{Dep. Matem\`atica Aplicada IV, EPSC,
     Universitat Polit\`ecnica de Catalunya, c/ Esteve Terradas 5, Castelldefels (Barcelona), Catalonia,
     Spain ({\tt comellas@ma4.upc.edu}).}
\author{Hongxiao Liu,  Zhongzhi Zhang}
\address{School of Computer Science and  Shanghai Key Lab of Intelligent Information Processing,
      Fudan  University, Shanghai 200433, China \\ ({\tt zhangzz@fudan.edu.cn}, {\tt 10210240077@fudan.edu.cn}).}


\begin{abstract}
In this paper we give an exact analytical expression for the number of spanning trees 
of an infinite family of  outerplanar, small-world and self-similar graphs.
This number is an important graph invariant
related to different topological and dynamic properties of the graph, such as
its reliability,  synchronization capability and diffusion properties.
The calculation of the number of spanning trees is a demanding and difficult task, in particular 
for large graphs, and thus there is much interest in obtaining closed expressions for relevant 
infinite graph families.
We have also calculated the spanning tree entropy of the graphs which we have
compared with those for graphs with the same average degree.
\end{abstract}



\end{frontmatter}

\section{Introduction }

Finding the number of spanning trees of graph is an old problem relevant 
in areas as diverse as mathematics,  chemistry, physics,  and
computer science.  This graph invariant is a  parameter 
related, for example,  to  the reliability of a network \cite{Co87}
its  synchronization \cite{TaMo06} and the study of random walks \cite{Ma00}. 
The number of spanning trees of any  finite graph can be computed
from the well known Kirchhoff's matrix-tree theorem as
the product of all nonzero eigenvalues of the
Laplacian matrix of the graph~\cite{GoRo01},  however
this is a demanding and difficult  task,  in particular for large graphs. 
Thus,  there has been much interest  in finding other methods to produce exact expressions 
for the number of spanning trees of relevant graph families such as grids
\cite{NiPa04}, lattices \cite{ShWu00} and Sierpinski
gaskets~\cite{ChChYa07,TeWa11}.

In this paper we give an exact analytical expression for the number of spanning trees 
of an  infinite family of outerplanar, small-world and self-similar graphs
which was introduced by two of the authors in~\cite{CoMi09}.
This family also has clustering zero and an exponential degree distribution.
The  combination  of modularity, small clustering coefficient,  small-world properties and
exponential degree distribution can be found in some networks associated with several social and technological 
systems~\cite{Ne03, FeJaSo01}.  Moreover, outer planarity is a relevant  property 
as it is known that many algorithms that are NP-complete 
for general graphs perform polynomial in outerplanar graphs~\cite{BrLeSp99}.


\section{The $M(t)$ graph family: construction and properties}

\begin{figure}[h]
\centerline{\includegraphics[width=1\linewidth]{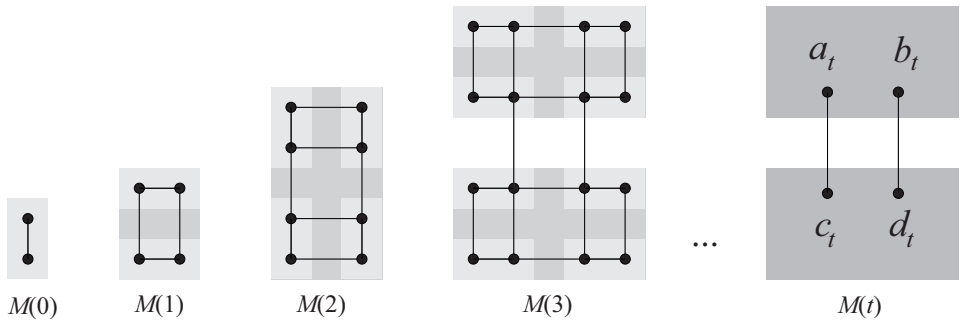}}
\caption{Graphs $M(t)$ produced at iterations $t=0,1, 2, 3$ and $t$.}\label{fig:recurrencia}
\end{figure}

This family of graphs  $M(t)$ was introduced in \cite{CoMi09} and can be
constructed following a recursive-modular method, see Fig.~\ref{fig:recurrencia}. 

For $t=0$, $M(0)$ has two vertices and one edge. 

For $t=1$, $M(1)$ is obtained from two graphs $M(0)$ connected by two new edges. 

For $t\ge 2$, $M(t)$ is obtained from two graphs $M(t\!-\!1)$ by connecting
them with two new edges. In each $M(t\!-\!1)$ the two vertices chosen are adjacent with maximum degree and
 have also been used  at  step $t\!-\!1$ to connect two $M(t\!-\!2)$.


 The number of vertices and edges of $M(t)=(V_t,E_t)$ are, respectively,
$|V_t|=2^{t+1}$  and  $|E_t|=3\cdot 2^{t}-2$.

The graph $M(t)$  is outerplanar, that
is,  the graph is planar and has an embedding in the plane such that 
all  vertices lie always in the exterior face of the graph while the edges never cross.
This family of graphs has the degrees of the endvertices positively
correlated --the graphs are assortative-- while
their  degree distribution  follows an exponential law
$P_{\rm cum}(k)= 2^{2-k}$.  See  \cite{CoMi09,CoMi09b} for more details and properties.

For large $t$, the diameter and  average distance of $M(t)$ scale logarithmically  with the order of the graph, 
meaning that the graph is small-world, although its clustering coefficient is zero.


\section{The number of spanning trees of $M(t)$}
In this section we make use of the construction method of $M(t)$ described in the previous section 
to obtain a set of recursive equations for the number of spanning trees and spanning subgraphs of 
the graph, which then can be solved by induction to find the exact number of spanning trees of $M(t)$ 
at any iteration step $t$.
For this calculation, we adapt the decimation technique given 
in~\cite{KnVa86}, which has been also considered to
find the number of spanning trees of the Sierpinski
gasket, the pseudofractal web, and several fractal lattices, 
see~\cite{ChChYa07,ZhLiWuZo10,ZhLiWuZo11}.

As shown in the previous section, $M(t)$ can be obtained from two $M(t\!-\!1)$ graphs 
interconnected by adding two new edges. 
As the four endvertices of these new edges are particularly relevant when counting the
number of spanning trees of $M(t)$,   we denote 
them as $a_t$, $b_t$,  $c_t$ and $d_t$ and we will refer to them as  {\em hubs}, 
see Fig.~\ref{fig:recurrencia}.

We recall here that a spanning subgraph of $M(t)$ is a subgraph with the same vertex set 
as $M(t)$  and a number of edges $|E'_t |$ such that $|E'_t | \leq  |E_t|$. 
A spanning tree of $M(t)$ is a spanning subgraph which is a tree and thus $|E'_t|= |V_t| - 1$.

Let $s(t)$ denote the total number of spanning trees of   $M(t)$ and
let $g(t)$ be the number of spanning subgraphs of $M(t)$ which consist of
two trees such that one of  the hub vertices used to form $M(t\!+\!1)$ belongs to one tree and the other 
is in the second tree (note that one of the trees might be just one hub).

To find  $s(t)$, we classify the different cases that contribute to $s(t)$ 
taking into account the connections between hub vertices of the two copies of $M(t\!-\!1)$ that
form $M(t)$.

We denote the existence of an edge between hubs $a_t$, $c_t$ by  $( a_t, c_t)$ 
and we represent that there is no edge between them by $\overline{( a_t, c_t)}$.
We use the same notation for  hubs $b_t$ and $d_t$. 

\begin{enumerate}
 \item If $\overline{(a_t,c_t)}$ and $\overline {(b_t, d_t)}$ then there is no contribution to $s(t)$.
\item If $(a_t,c_t)$ and $\overline {(b_t, d_t)}$, the number of spanning trees that contribute to $s(t)$ is  
\begin{equation}\label{equation:s1}
s(t\!\! -\!\!1)^2.
\end{equation}
\item If $\overline {(a_t,c_t)}$ and $(b_t, d_t)$, like in the previous case,  the number of spanning trees is  
\begin{equation}\label{equation:s2}
s(t\!\! -\!\!1)^2.
\end{equation}
\item If $(a_t,c_t)$ and $(b_t, d_t)$.
 In this case the number of spanning trees contributing to $s(t)$ is 
 $s(t\!\! -\!\!1)g(t\!\! -\!\!1)$, i.e. the product of the number of spanning trees  
 of the copy $M(t\!-\!1)$  containing hubs $a_t$ and $b_t$ with 
$g(t\!\! -\!\!1)$, which counts  the number of spanning subgraphs of  $M(t\!-\!1)$ such that each contains the  
 hubs $c_t$ and $d_t$ but does not have a path connecting them because they belong to different trees 
 of the spanning graph.
We have  also to consider the symmetrical case with respect to the sets $\{a_t,b_t\}$ and $\{c_t,d_t\}$.
Thus, the  total number of spanning trees that corresponds to this case is
\begin{equation}\label{equation:s3}
 2s(t\!\! -\!\!1)g(t\!\! -\!\!1).
\end{equation}
\end{enumerate}

By adding (\ref{equation:s1}), (\ref{equation:s2}) and (\ref{equation:s3}) we find that the total number of 
spanning trees of $M(t)$ is
\begin{equation}\label{equation:sst}
s(t)= 2s(t\!\! -\!\!1)^2+  2s(t\!\! -\!\!1)g(t\!\! -\!\!1)
\end{equation}
Now we need to obtain a recursive relation for $g(t)$.
A similar reasoning than above leads to the following equation
\begin{equation}\label{equation:gt}
 g(t)= s(t\!\! -\!\!1)^2+ 2s(t\!\! -\!\!1)g(t\!\! -\!\!1).
\end{equation}

Subtracting equation (\ref{equation:gt}) from equation (\ref{equation:sst}) we obtain 
\begin{equation*}
s(t)- g(t) = s(t\!\! -\!\!1)^2
\end{equation*}
and thus 
\begin{equation}\label{equation:gtt}
 g(t\!\! -\!\!1) = s(t\!\! -\!\!1)-s(t\!\! -\!\!2)^2.
\end{equation}
Introducing (\ref{equation:gtt}) into (\ref{equation:sst}) we have:
\begin{equation*}
s(t)= 4s(t\!\! -\!\!1)^2-2s(t\!\! -\!\!1)s(t\!\! -\!\!2)^2 \quad \text{with}\quad s(0)=1\quad \text{and}\quad s(1)=4.
\end{equation*}
To solve this equation, we rewrite it as
\begin{equation}\label{equation:st}
\frac{s(t)}{s(t\!\! -\!\!1)^2}= 4-2\frac{1}{\frac{s(t\!\! -\!\!1)}{s(t\!\! -\!\!2)^2}}
\end{equation}
and, if we introduce 
\begin{eqnarray}\label{equation:qtdef}
q(t)&=&s(t) / s(t\!\! -\!\!1)^2, \quad t\ge 2,\\
q(1)&=&4,\nonumber
\end{eqnarray}
equation (\ref{equation:st}) becomes
\begin{equation}
 q(t) = 4 - \frac{2}{q(t\!\! -\!\!1)}, \quad  \text{with}\quad t\geq 2,
\end{equation}
which can be solved to obtain
\begin{equation}\label{equation:qt}
 q(t)=2-\sqrt{2} + \frac{2 \sqrt{2}} {1 -(3-2\sqrt{2})^t}. 
\end{equation}
From equation~ (\ref{equation:qtdef}) we have:
\begin{eqnarray*}
s(t)\!\!\!\!&=&\!\!\!\!q(t) s(t\!\! -\!\!1)^2=q(t) q(t\!\! -\!\!1)^2 s(t\!\! -\!\!2)^4 = \cdots  \\
      \!\! \! \!&=&\!\!\!\!q(t)  q(t\!\! -\!\!1)^2 q(t\!\! -\!\!2)^{2^2} q(t\!\! -\!\!3)^{2^3} \cdots q(t\!\! -\!\!i)^{2^i}  
        \cdots q(1)^{2^{t-1}},
\end{eqnarray*}
and thus, using~(\ref{equation:qt}), we can write:

\begin{theorem}\label{theo:main}
The number of spanning trees of $M(t)$, $t \geq 1$, is
 \begin{equation*}
 s(t) = \prod_{i=1}^{t}\left[ 2-\sqrt{2} + \frac{2 \sqrt{2}} {1 -(3-2\sqrt{2})^i}\right]^{2^{t-i}}.
 \end{equation*}

\end{theorem}

\section{Spanning tree entropy of $M(t)$}
After having an exact  expression for the number of spanning trees
of $M(t)$, we can calculate its spanning tree entropy, which
is defined as in~\cite{Ly05,Wu77}:
\begin{equation*}
\label{eq:en} h=\lim_{V_t \to \infty}\frac{\ln s(t)}{V_t}.
\end{equation*}
and we obtain
$$h\simeq 0.657.$$

We can now compare this asymptotic value of the entropy of the spanning trees of $M(t)$
with those of other graph families with the same average degree:
As an example, the value for the honeycomb lattice is 0.807~\cite{Wu77} and the
4-8-8 (bathroom tile) and 3-12-12 lattices have  entropy values 0.787 and 0.721,
respectively~\cite{ShWu00} while the spanning tree entropy for Hanoi graphs is 
0.677~\cite{ZhLiWuCo12} .
Thus, the asymptotic value of the entropy of the spanning trees of  $M(t)$ is the lowest known for
graphs with average degree 3.
This means that the number of spanning trees in $M(t)$,
although growing exponentially, do it at a
lower rate than graphs with the same average degree.
This result  suggests that the family of graphs $M(t)$ would be less reliable
to a random removal of edges  than  the graphs mentioned above (as they have more spanning trees).
However, the degree distribution also affects the reliability of a graph and
it is well known that graphs with a scale-free distribution of degrees are much more 
resilient than homogeneous graphs (like regular and exponential graphs, where all vertices are  
statistically identical), see~\cite{Tu00}. 
Thus, the  large inhomogeneity in the degree distribution of the graphs $M(t)$ 
might increase their robustness and it would be of interest to study, for these and other relevant graph families, 
the connections among spanning tree entropy, degree distribution and other topological parameters 
like degree correlation, in relation to their dynamic properties (reliability, synchronizability, etc.)

\section{Conclusion}
In this paper we find the number of spanning trees 
of the $M(t)$ graph family by using a method, based on their self-similar structure,
which allows us to obtain an analytical exact expression valid for any number of iterations.
From  the  number of spanning trees  we find
their spanning tree entropy which is lower than in other graph families  with the same average degree.

\section*{Acknowledgments}

Z.Z. and H.L. are supported by the National Natural Science
Foundation of China under grant No. 61074119. F.C. and A.M. are supported by the Ministerio de Economia y Competitividad, Spain, and the European Regional Development Fund under project MTM2011-28800-C02-01 and partially supported by the Catalan Research Council under grant 2009SGR1387.


\end{document}